\documentclass{amsart}%
\usepackage{amsfonts}
\usepackage{amsmath}
\usepackage{amssymb}
\usepackage{graphicx}%
\providecommand{\U}[1]{\protect\rule{.1in}{.1in}}
\newtheorem{theorem}{Theorem}
\theoremstyle{plain}

\newtheorem{corollary}{Corollary}

\newtheorem{definition}{Definition}

\newtheorem{lemma}{Lemma}

\newtheorem{proposition}{Proposition}
\newtheorem{remark}{Remark}

\numberwithin{equation}{section}
\begin{document}
\title[Cauchy problem for NLKG in modulation spaces with noninteger powers]
{Cauchy problem for NLKG in modulation spaces with noninteger powers}
\author{Huang Qiang}
\address[Huang Qiang]{Department~of~Mathematics,~Zhejiang~University,~Hangzhou~310027,~PR~China}
\email[Huang Qiang]{huangqiang0704@163.com}
%\urladdr{http://www.authorone.oneuniv.edu}
\author{Fan~Dashan}
\address[Fan~Dashan]{Department~of~Mathematics,~University~of~Wisconsin-Milwaukee,~Milwaukee,~WI~53201,~USA}
\email[Fan~Dashan]{fan@uwm.edu}
%\urladdr{http://www.authortwo.twouniv.edu}
\author{Chen~Jiecheng}
\address[Chen~Jiecheng]{Department~of~Mathematics,~Zhejiang~Normal~University,~Jinhua~321004,~PR~China}
\email{jcchen@zjnu.edu.cn}
%\urladdr{http://www.authorthree.threeuniv.edu}
%\thanks{Thanks for Author One.}
%\thanks{Thanks for Author Two.}
%\thanks{This paper is in final form and no version of it will be submitted for publication elsewhere.}
\date{November 11, 2014}
\subjclass[2000]{35A01, 35A02, 42B37}
\keywords{Modulation spaces, Nonlinear Klein-Gordon equation, Cauchy problem, noninteger power}
\dedicatory{ }
\begin{abstract}
In this paper, we consider the Cauchy problem for the nonlinear Klein-Gordon
equation whose nonlinearity is $|u|^{k}u$ in the modulation space, where $k$
is not an integer. Our method can be applied to other equations whose
nonlinear parts have regularity estimates. We also study the global solution
with small initial value for the Klein-Gordon-Hartree equation. By this we
can show some advantages of modulation spaces both in high and low
regularity cases.

\end{abstract}
\maketitle

\section{Introduction and main results}

\hspace{6mm} In this paper, we study the Cauchy problem for the following
nonlinear Klein-Gordon equation (NLKG):
$$
u_{tt}+(I-\triangle )u=\pm |u|^{k}u,~~~~~~~~u(0)=u_{0},u_{t}(0)=u_{1},\eqno%
(1.1)
$$%
where
\[
k\in(0,+\infty)\backslash\mathbb{Z^{+}},u_{tt}=\partial ^{2}/\partial ^{2}t,
\]
and $\ $%
\[
\triangle =\partial ^{2}/\partial ^{2}x_{1}+...+\partial ^{2}/\partial
^{2}x_{n}
\]
is the Laplace operator. It is well known that the NLKG has the following
equivalent integral form:
$$
u(t)=K^{\prime }(t)u_{0}+K(t)u_{1}-\int_{0}^{t}K(t-\tau )|u|^{k}ud\tau ,\eqno%
(1.2)
$$%
where we denote $\ \omega =(I-\Delta )$ $\ $and
\[
K(t)=\frac{\sin t\omega ^{\frac{1}{2}}}{\omega ^{\frac{1}{2}}}%
,~~~~~~K^{\prime }(t)=\cos t \omega ^{\frac{1}{2}}.
\]

\bigskip

The aim of this paper is to study the local and global well posedness of
NLKG in modulation spaces. The modulation space \ $M_{p,q}^{s}$ \ was
originally introduced by Feichtinger in \cite{F} , where its definition is
based on the short-time Fourier transform and the window function.
Feichtinger's initial motivation was to use the modulation space to measure
smoothness for some function or distribution spaces. Since then, this space
was received an extensive study on its analysis/topological constructure and
algebraic properties. See, for example, \cite{GS,M2,M3,M1} and the
references therein for more details. Later, people found that this space is
a good working frame to study certain Cauchy problems of nonlinear partial
differential equations. To this end, Wang and Hudzik gave another equivalent
definition in \cite{BHH} \ using the frequency-uniform-decomposition
operators. With this discrete definition, they was able to consider the
global solutions for nonlinear Schr\"{o}dinger equation and nonlinear
Klein-Gordon equation in the space $M_{p,q}^{s}.$ Following their pioneer
work, we can find a lot of research papers in the literature that address
various harmonic analysis and PDE problems on the modulation spaces. In the
following, we list a few of these results, among numerious of papers. Gr\"{o}%
bner in his PH.D. thesis \cite{G} introduced the \ $\alpha $ \ modulation
spaces that reveals the essential connection between the modulation spaces
and the Besov spaces; Han and Wang\cite{HW} followed Gr\"{o}bner's idea to
give a discrete version of the $\ \alpha $ \ modulation space based on the
frequency-uniform-decomposition, and they obtained more properties related
to this space; Feichtinger, Huang and Wang \cite{FHW} studied the trace
operator in modulation, $\alpha $-modulation and Besov spaces. Also, for PDE
problems, Wang, Zhao, and Guo \cite{BZG} studied the local solution for
nonlinear Schr\"{o}dinger equation and Navier-Stokes equations; Wang and
Huang \cite{BH} obtained the local and global solutions for generalized KdV
equations, Benjamin-Ono and Schr\"{o}dinger equations; Tsukasa Iwabuchi
studied the local and global solutions for Navier-Stokes equations, as well
as the heat equations (see \cite{TI}). However, we observe that the
nonlinear parts of above mentioned equations is either $\ |u|^{k}u$ \ with $%
k\in \mathbb{Z^{+}}$ \ or a multi-linear function $F(u,..,u)$. The reason of
such a restriction is that all the estimates in \ \cite{BZG}, \cite{BH},
\cite{TI} are based on the algebra property
\[
\left\Vert u^{2}\right\Vert _{M_{p,1}^{s}}\preceq \left\Vert u\right\Vert
_{M_{p,1}^{s}}^{2},
\]%
that causes that the exponent \ $k$ \ must be a positive integer in the
nonlinear term \ $|u|^{k}u.$ In a recent article \cite{RMW}, Ruzhansky,
Sugimoto and Wang stated some new progresses on the modulation spaces. In
the same article they posed three open questions. One of these questions is
to study nonlinear PDE whose nonlinear term \ $|u|^{k}u$ \ has a non-integer
$k\in (0,+\infty )$.

Motivated by the above question, in this paper, we will give a partial
answer to the above problem. Before we attack the problem, let us briefly
describe some obvious difficulties in handling this problem. Unlike the
Lebesgue spaces or the Besov spaces, we do not have the Littlewood-Paley
theory or any of its analog on the modulation spaces. So we can only use its
algebra property and local analysis on its windows in the frequency spaces.
However, both of these two tools handle only the case of integer $\ k$.
These bring the main difficulty in our problem. To achieve our aim,
fortunately we observe that modulation spaces and Besov spaces can be
embedded each other. So, we can use this property to deduce the problem on
the modulation spaces to those on the Besov spaces and then transfer the
obtained estimates back to the modulation spaces. Of course, during this
embedding process, we might loss some regularity. So this method might be
applied in some equations whose nonlinear parts have regularity estimates.
These equations include the Klein-Gordon equation, the heat equation and
some other equations with regularity estimates.

Now we state our main theorems in the following.

\begin{theorem}
Let $s\in \mathbb{R},k\geq[s],2<p<\infty ,$and $1\leq q<\infty .$ Assume that \ $%
q $ \ satisfies the following conditions:
$$
max\{1+\frac{n}{2}-\frac{n}{q},\frac{n}{2}-\frac{1}{k}[1-(\frac{1}{q}-\frac{1%
}{2})n]\}<s<\frac{n}{2}\eqno(1.3)
$$%
when $q<2$ $\ $and $n(\frac{1}{q}-\frac{1}{2})<1$;
\end{theorem}

\textit{and}
$$
max\{1,\frac{n}{q^{\prime }}-\frac{1}{k}[1-n(\frac{1}{2}-\frac{1}{q})]\}<s<%
\frac{n}{q^{\prime }}\eqno(1.4)
$$%
\textit{when }$q>2,n(\frac{1}{2}-\frac{1}{q})<1$\textit{. }

\textit{For any }$(u_{0},u_{1})\in M_{2,q}^{s}\times ,M_{2,q}^{s-1}$\textit{%
\ , there exists a }$T>0$\textit{\ such that the equation (1.1) has a
unique solution in the space }%
$$
L^{\rho }(0,T;M_{p,q}^{s-\beta })\bigcap L^{\infty }(0,T;M_{2,q}^{s}),\eqno%
(1.5)
$$%
\textit{where }$\rho $\textit{\ and }$\beta \ \ $\textit{are any real
numbers satisfying }$\ \rho \beta =\frac{n+1}{n-1}$\textit{\ and }$\beta \in
(0,\frac{n+1}{2n-2}).$

Since $\beta >0\ \ $is arbitrary, $L^{\rho }(0,T;M_{p,q}^{s-\beta })$ might
be closed to the space $L^{\infty }(0,T;M_{p,q}^{s}),$ which has more
regularity and better integrality .

\begin{remark}
When consider local well-posedness of NLKG in Besov space $B_{p,q}^{s}$, the
domain of $p,q$ \ is $\ 1\leq q\leq \infty $ \ and $(\frac{1}{2}-\frac{1}{p}%
)n\in \lbrack 0,1)$. Comparing this to Theorem 1, we can see that the
domain of $\ p$ \ in the modulation space is similar to the domain of $q$ \
in the Besov space, and $q$ behaves in modulation spaces quite like the
performance of $\ p$ \ in the Besov space. One can easily understand this
nature, when one compares the embedding relations between the modulation
spaces (see 2.1) to the Sobolev embedding in Besov spaces.
\end{remark}

In Theorem 1, we use an auxiliary space $L^{\infty }(0,T;M_{2,q}^{s})$.
The purpose is that we want to expand the domain of $p,q$ \ by the Stricharz
estimate. Actually, we can also prove the unconditional local well posedness
in $L^{\gamma }(0,T;M_{p,q}^{s})$ without the auxiliary space. As we all
know, the unimodular semigroup \ $e^{it\left\vert \Delta \right\vert
^{\alpha /2}}$ is bounded on the $\ L^{p}$ space or on the Besov space $%
B_{p,q}^{s}$ \ if and only if $p=2$ (or \ $\alpha =1$ at \ $n=1).$ So if we
want to obtain the local well posedness in $\ B_{p,q}^{s}$ when $p\neq 2$,
we should use the Sobolev space \ $H^{s}$ as an auxiliary space. On the
other hand, for the Sobolev space $H^{s}$ alone, it is difficult to estimate
nonlinear part when $s<\frac{n}{2}$, so it also needs the Besov space as an
auxiliary space. But in modulation spaces, we are able to obtain the $%
M_{p,q}^{s}-$boundedness for the unimodular semigroup and we can also
estimate the nonlinear part in the low regularity case. As these advantages,
the next corollary shows that we can obtain the following unconditional
local well posedness in modulation spaces. This feature is not available
either on the Besov spaces, or on the Sobolev spaces, when one studies the
problem in the low regularity case.

\begin{corollary}
Let $1\leq q<\infty ,2\leq p<\infty $,$\ s\in \mathbb{R}$  and $k\geq[s]$ . Assume that
they satisfy $q\in \lbrack p^{\prime },p],(1-\frac{2}{p})n<1$ and
$$
\max \{1-n(\frac{1}{q}-\frac{1}{p}),\frac{n}{q^{\prime }}-\frac{1}{k}\}<s<%
\frac{n}{q^{\prime }}-\frac{1}{k}(1-\frac{2}{p})n.\eqno(1.6)
$$%
Then for any initial value $(u_{0},u_{1})\in M_{p,q}^{s}\times
,M_{p,q}^{s-1} $, there exists a $T>0$ such that \textit{the equation (1.1)}
has an unique solution in $L^{\gamma }(0,T;M_{p,q}^{s})$ for any $\gamma
\geq k+1$.
\end{corollary}

\begin{remark}
Throughout the proof, we will use the Stricharz estimate and a nonlinear
estimate in Besov spaces. Both of these two tools requires $2<p<\infty $.
So, for the case $1\leq p\leq 2$, we are not able to obtain the
well-posedness in the Besov spaces. However, such a restriction can be
removed when we study the same problem on the modulation spaces. In the next
theorem, we indeed obtain the solution in $M_{p,q}^{s}$ when $1\leq p\leq 2$%
. Recall that in the Besov spaces $\ B_{p,q}^{s}$\ , if we want to use
Sobolev embedding to control the norm \ $B_{p,q}^{s}$\ \ by \ the norm \ $%
B_{p_{1},q}^{s_{1}}$ \ with \ $p_{1}<p$, it needs more regularity $s_{1}>s$.
But in modulation spaces $M_{p,q}^{s}$, we have uniformly estimate for the
index $p$ (see (2.1)) which has no influence to the regularity index $s$.
Hence we can use this embedding property and boundedness of unimodular
semigroup to solve the problem in the case $1\leq p\leq 2$. Specifically, we
have the following theorem:
\end{remark}

\begin{theorem}
Let $1<p\leq 2$, $1\leq q<\infty $ $\ $$\ s\in \mathbb{R}$  and $k\geq[s]$.  Assume that they
satisfy that $\ q\in \lbrack p,p^{\prime }],(\frac{2}{p}-1)n<1$, and
$$
\max \{1-n(\frac{1}{q}-\frac{1}{p^{\prime }}),\frac{n}{q^{\prime }}-\frac{1}{%
k}\}<s<\frac{n}{q^{\prime }}-\frac{1}{k}(1-\frac{2}{p^{\prime }})n.\eqno%
(1.7)
$$%
Then the same conclusion as Corollary 1 holds .
\end{theorem}

\begin{remark}
In the case $p=2$, the condition $q\in \lbrack p^{\prime },p]$ in the above
theorem means $q=2$. In \cite{BHH}, Wang and Hudzuk proved that the space $\
M_{2,2}^{s}$ \ is equivalent to the Sobolev space $H^{s}$, so this result is
not interesting. Actually, when $p=2$, it is not necessary to choose $q=2$.
The range of $q$ \ for $p=2$ can be wider. We will address this special case
in Remark 6.
\end{remark}

Now we turn to state the global solution of NLKG in modulation spaces.

\begin{theorem}
Let $1\leq q<\infty ,2\leq p<\infty $ $\ s\in \mathbb{R}$  and $k\geq[s]$ . Assume that they satisfy $q\in \lbrack p^{\prime },p],(1-%
\frac{2}{p})n<1-\alpha $ and
$$
\max \{1-\frac{\alpha }{2}-n(\frac{1}{q}-\frac{1}{p}),\frac{n}{q^{\prime }}-%
\frac{1}{k}(1-\alpha )+\frac{\alpha }{2}\}<s<\frac{n}{q^{\prime }}-\frac{1}{k%
}(1-\frac{2}{p})n+\frac{\alpha }{2}\eqno(1.8)
$$%
where
$$
\alpha =\theta (n+1)(\frac{1}{2}-\frac{1}{p}),~~~~~~~~\delta =\theta (n-1)(%
\frac{1}{2}-\frac{1}{p}).\eqno(1.9)
$$%
for $\ $some $\ \theta \in \lbrack 0,1]$. \ In addition, assume $q\in
\lbrack \gamma ^{\prime }.\gamma ],\gamma \geq \frac{2}{\delta }$, and $k>%
\frac{4}{\theta }+\frac{2}{n}$. Then there exists a small $\nu >0$ such that
for any $\Vert u_{0}\Vert _{M_{p,q}^{s}}+\Vert u_{1}\Vert
_{M_{p,q}^{s-1}}\leq \nu $, equation (1.1) has an unique global solution
$$
u\in L^{\infty }(R;M_{p,q}^{s})\bigcap L^{k+2}(R;M_{p,q}^{s-\frac{\alpha }{2}%
}).\eqno(1.10)
$$
\end{theorem}

\begin{remark}
When we choose $\theta =1$ in Theorem 3, we can find the global solution
with small initial value for $k>4+\frac{2}{n}$. In \cite{BHH}, Wang and
Hudzik proved that if $k$ is an integer, the global existence interval is $%
k\geq \frac{4}{n}$ which is wider then ours. That is because during the
embedding between Modulation spaces and Besov spaces we lost some
regularity. Hence we need more power of $u$ to guarantee more regularity to
make up this lost.
\end{remark}

By far, we find the local solution for $1<p<\infty $ and the global solution
for $2<p<\infty $ when $k$ is not integer. In all cases, nonlinear estimate
in modulation spaces relies heavily on the corresponding nonlinear estimate
in Besov spaces and some regularity is lost. However, we obtain two
advantages in the modulation spaces and they are quite unique comparing to
the results on the Besov space. First, in Corollary 1, we obtain the
unconditional local well posedness with low regularity. Second. in Theorem
2, we solve the problem in the case $1<p<2$. Moreover, if the nonlinearity
is a multi-linear function, we can find that the modulation spaces have many
other advantages. Below, we will use the nonlinear Klein-Gordon-Hartree
equation(NLKGH):
$$
u_{tt}+(I-\triangle )u+(|x|^{-\mu }\ast u^{2})u=0,u(0)=u_{0},u_{t}(0)=u_{1}%
\eqno(1.11)
$$%
to illustrate these advantages.

We consider the global solution with small initial value for equation (1.11)
in $\mathbb{R}\times \mathbb{R}^{3}$, and compare the result with the same
solution in Besov spaces obtained in \cite{MZ}. As we mentioned before, the
role of the index $q$ \ in the modulation space is significant. The
regularity index $s$ might depend on \ $q$ \ in some estimates, for instance
see (2.2). In order to obtain a good time-space estimate, in \cite{BHH} Wang
and Hudzik gave up a traditional method of dual estimate. They introduced
the space $l_{\Box }^{s,q}(L^{r}(0,T;L^{p}(R^{n})))$ (see Definition 1 )
to replace the standard modulation space. Then, for $q\in \lbrack \gamma
^{\prime }.\gamma ]$, they were able to invoke the Minkowski inequality to
obtain some time-space estimates in $L^{r}(0,T;M_{p,q}^{s})$. Look back to
Theorem 3, in the proof we need the embedding between modulation spaces
and Besov spaces to solve the case $k\notin \mathbb{Z}$, so we work on the
space $L^{r}(0,T;M_{p,q}^{s})$ whose Stricharz estimate needs to restrict to
$q\in \lbrack \gamma ^{\prime }.\gamma ]$. But for the equation (1.11), we
want to use the relation between $s$ and $q$ to extend the domain of $\mu $.
So we hope that the restrictions on $q$ \ are as less as impossible. To
achieve this target, we choose the space $l_{\Box
}^{s,q}(L^{r}(R;L^{p}(R^{n})))$ to find the global solution and establish
the following result.

\begin{theorem}
Suppose that $(u_{0},u_{1})\in M_{2,q}^{s}\times M_{2,q}^{s-1},$ where $%
1<q<\infty ,s\geq 0$, $\alpha $ and $\delta $ are defined in (1.9). Assume
that the domain of $p$ satisfies
\[
(1-\frac{2}{p})\in \lbrack \frac{1}{2\theta (n-1)},\frac{1}{\theta (n-1)}).
\]%
For $\ 2n(1-\frac{2}{p})\leq \mu \leq 2(s+\frac{n}{q})+1-2\alpha -n$, we can
find a constant $\varepsilon >0$ for which if
\[
\Vert u_{0}\Vert _{M_{2,q}^{s}}+\Vert u_{1}\Vert _{M_{2,q}^{s-1}}\leq
\varepsilon
\]%
then, for equation (1.3), there exists a unique global
$$
u\in l_{\Box }^{s,q}(L^{\infty }(R;L^{p}(R^{n})))\bigcap l_{\Box }^{s-\frac{%
\alpha }{2},q}(L^{4}(R;L^{p}(R^{n})).\eqno(1.12)
$$
\end{theorem}

\begin{remark}
In \cite{MZ}, Miao and Zhang studied equation (1.11) on the Besov spaces. In
the case \ $n=3,$ they showed that the exponent $\mu $ must satisfy
\[
\mu =\frac{6(s+1)}{3+\eta }
\]%
and
\[
\mu \geq \frac{6}{2+\eta }
\]%
where $\eta \in \lbrack 0,1]$. \ These requirements imply $s\geq \frac{1}{%
2+\eta }$. So the minimum regularity should be $s=\frac{1}{3}$ when we
choose $\eta =1$. But on the modulation space in Theorem 1.9, if we let $%
\theta =1$ , we can choose $p$ such that $2n(1-\frac{2}{p})=\frac{3}{2}$
when \ $n=3$. Then, if we choose $q$ closed to 1, it is not difficult to
find that the minimum value of $\ s\ $can be $0.$ We observe that in \cite%
{BHH}, Wang and Hudzik proved that $M_{p,q}$ has no derivative regularity
for any $0<p,q\leq \infty $. Hence, our result is another form of low
regularity for global solution of NLKGH. Second, If one wants to obtain some
high regularity estimate for this equation, in \cite{MZ} the domain of $\mu $
is $\frac{3}{2}(s+1)\leq \mu \leq 2(s+1)$ when the authors take the Besov
space as the working space. Checking the domain of \ $\mu $ on modulation
spaces in Theorem 4, clearly it is larger, since the low bound is fixed
which is independent on $s$. So, both in low regularity and high regularity
cases, modulation spaces seems better than Besov spaces.
\end{remark}

We are not surprising that the modulation space have these advantages
comparing to the Besov space. In the Besov spaces, many estimates, such as
the admissible pairs, H\"{o}lder's inequality, boundedness of fractional
integral operator, the Sobolev embedding, etc., rely all on the exponent $p$%
, while the index $q$ is dummy. But in the modulation spaces, the admissible
pair relies on $p$, the Sobolev embedding relies on $q$ \ (see (2.2)), H\"{o}%
lder's inequality and boundedness of fractional integral operator rely on
both $\ p$ and $q$. Moreover, we have the uniformly estimate for the index $%
p $ (see (2.1)). In other words, working in the Besov spaces, one needs to
give too much restrictions of $p$, and $q$ plays no role. In the modulation
spaces, $p$ and $q$ share these restrictions together, and both $p$ and $q$
have uniformly estimates.

The proofs of theorems will be represented in the third section.

\section{Preliminaries}

\hspace{6mm}In this section we recall the definitions and some properties of
the modulation space and Besov space. Also, we will prove several lemmas,
particularly a key lemma to estimate $|u|^{k}u$ in the modulation space when
$k$ is not an integer.

\begin{definition}
(Modulation spaces) Let $\{\varphi _{k}\}\subset C_{0}^{\infty }(R^{n})$ be
a partition of the unity satisfying the following conditions:
\[
supp\varphi \subset \{\xi \in R^{n}|\mid \xi \mid \leq \sqrt{n}\},\sum_{k\in
Z^{n}}\varphi _{k}(\xi )=1,
\]%
for any $\xi \in R^{n}$, where $\varphi _{k}(\xi ):=\varphi (\xi -k)$. For
each \ $k\in
%TCIMACRO{\U{2124} }%
%BeginExpansion
\mathbb{Z}
%EndExpansion
^{n},$ denote a local square projection \ $\Box _{k}$ \ on the frequency
space by \
\[
\Box _{k}:=\mathcal{F}^{-1}\varphi _{k}\mathcal{F},
\]%
where $\ \mathcal{F}$ \ and \ $\mathcal{F}^{-1}$denote the Fourier transform
and its inverse, respectively. By this frequency-uniform decomposition
operator, we define two kinds of modulation spaces, for $\ 0<p,q\leq \infty $
and $\ s\in
%TCIMACRO{\U{211d} }%
%BeginExpansion
\mathbb{R}
%EndExpansion
,$ by
\[
M_{p,q}^{s}(R^{n}):=\{f\in S^{\prime }(R^{n}):\Vert f\Vert
_{M_{p,q}^{s}(R^{n})}=(\sum_{k\in Z^{n}}\langle k\rangle^{sq}\Vert \Box _{k}f\Vert
_{p}^{q})^{\frac{1}{q}}<\infty \}
\]%
and
\[
l_{\Box }^{s,q}(L^{r}(0,T;L^{p}(R^{n}))):=\{f(t,\cdot )\in S^{\prime
}(R^{n}):\Vert f\Vert _{l_{\Box }^{s,q}(L^{r}(0,T;L^{p}(R^{n}))})<\infty \},
\]%
where \
\[
\Vert f\Vert _{l_{\Box }^{s,q}(L^{r}(0,T;L^{p}(R^{n})))}=(\sum_{k\in
Z^{n}}\langle k\rangle^{sq}\Vert \Box _{k}f\Vert _{L^{r}(0,T;L^{p}(R^{n}))}^{q})^{\frac{1%
}{q}}
\]%
and $\langle k\rangle:=(1+|k|^{2})^{\frac{1}{2}}$(See \cite{BHH} for details). If the
domain of $t$ is $(-\infty ,+\infty )$, we donate $l_{\Box
}^{s,q}(L^{r}L^{p})$ for convenience. The space $l_{\Box
}^{s,q}(L^{r}(0,T;L^{p}(R^{n})))$ was first introduced by Planchon \cite{P1},%
\cite{P2} when he studied the nonlinear Schr\"{o}dinger equation and the
nonlinear wave equation. In the definition, the order of $L^{r}$ norm and $%
l^{q}$ norm is changed. This change seems important in modulation spaces. As
we know, $q$ is a very important index in modulation spaces which can impact
the regularity. So, in many cases, we should deal with $q$ carefully and
choose $l^{q}$ norm in the last step. Moreover, we will recall some
properties of modulation spaces which will be useful in this paper. More
details can be found in \cite{BHH}.
\end{definition}

In the following content, if no special explanation, we always assume that
\[
s,s_{i}\in R,1\leq p,p_{i},q,q_{i}\leq \infty .
\]

\begin{proposition}
(Isomorphism \cite{BHH}).\newline
Let $\ 0<p,q\leq \infty ,s,\sigma \in R$. For the Bessel potential $%
J_{\sigma }=(I-\triangle )^{\frac{\sigma }{2}},$ the mappings
\[
J_{\sigma }:M_{p,q}^{s}\rightarrow M_{p,q}^{s-\sigma }~~
\]%
and
\[
J_{\sigma }:~l_{\Box }^{s,q}(L^{r}(0,T;L^{p}(R^{n})))\rightarrow l_{\Box
}^{s-\sigma ,q}(L^{r}(0,T;L^{p}(R^{n})))
\]%
are isomorphic mappings.
\end{proposition}

\begin{proposition}
(Embedding, \cite{BHH}).\newline
$M_{p_{1},q_{1}}^{s_{1}}\subset M_{p_{2},q_{2}}^{s_{2}}$ ~~and ~~$l_{\Box
}^{s_{1},q_{1}}(L^{r}(0,T;L^{p_{1}}))\subset l_{\Box
}^{s_{2},q_{2}}(L^{r}(0,T;L^{p_{2}}))$ if .\\

(i)~\text{\ }$s_{1}\geq s_{2},~0<p\leq p_{2},~0<q_{1}\leq
q_{2}$\qquad\qquad\qquad\qquad\qquad\qquad\qquad\qquad\qquad(2.1)\\

(ii)~$q_{1}>q_{2},~s_{1}>s_{2},~s_{1}-s_{2}>n/q_{2}-n/q_{1}$\qquad\qquad\qquad\qquad\qquad\qquad\qquad\qquad(2.2)

\end{proposition}

In (2.1), we can see that both $p$ and $q$ have uniform estimates, and from
(2.2) we can find that the condition on \ $q$ \ is similar to the Sobolev
embedding. Now we need the relationship between modulation spaces and Besov
spaces.

\begin{lemma}
(Embedding with Besov spaces, \cite{BHH})\newline
Assume $B_{p,q}^{s}$ is the Besov spaces, and $1\leq p,q\leq \infty ,s\in
\mathbb{R}$. We have the following embedding:
$$
M_{p,q}^{s+\sigma (p,q)}\subset B_{p,q}^{s},~~~~\sigma (p,q)=max(0,n(\frac{1%
}{p\wedge p^{\prime }}-\frac{1}{q}))\eqno(2.3)
$$%
$$
B_{p,q}^{s+\tau (p,q)}\subset M_{p,q}^{s},~~~~\tau (p,q)=max(0,n(\frac{1}{q}-%
\frac{1}{p\vee p^{\prime }})).\eqno(2.4)
$$
\end{lemma}

Moreover, since we will embed the modulation spaces into the Besov spaces,
we need employ some nonlinear estimates in Besov spaces, particularly the
estimate on the nonlinear term $u^{k}.$ Recall that such estimate in Besov
spaces \ $B_{p,q}^{s}$ \ has a long history. Cazenave obtained the case $%
0<s<1$ in \cite{C}, and Cazenave and Weissler obtained the case $1<s<\frac{N%
}{2}$ in \cite{CW}. Later, Wang proved a general case in \cite{B}. Our proof
will be based on Wang's result. Since all results on the Besov space \ $%
B_{p,q}^{s}$ \ are stated for the case $q=2$, we need the following
embedding to obtain information for all $1\leq q<\infty $, as we will handle
all \ $q$ \ in the space \ $M_{p,q}^{s}.$

\begin{proposition}
(\cite{T}) Let $\epsilon>0$, for any $1\leq p,q_{1},q_{2}\leq\infty$, then
we have
$$
B_{p,q_{1}}^{s+\epsilon}\subset B_{p,q_{2}}^{s}\eqno(2.5)
$$
\end{proposition}

\begin{lemma}
(Nonlinear estimate in Besov space)\newline
Suppose $2\leq p<\infty ,1\leq q\leq \infty $ \ and $0\leq \delta
<s<s_{1}<\infty $,$[s-\delta]\leq k$ If they satisfy
$$
k(\frac{1}{p}-\frac{s}{n})+\frac{1}{p}-\frac{\delta }{n}=\frac{1}{p^{\prime }%
},~~~~\frac{1}{p}-\frac{s}{n}>0\eqno(2.6)
$$%
then we have
$$
\Vert |u|^{k}u\Vert _{B_{p^{\prime },q}^{s-\delta }}\preceq \Vert u\Vert
_{B_{p,q}^{s_{1}}}^{k+1}\eqno(2.7)
$$
\end{lemma}

\textbf{Proof: } When $q=2$, in \cite{BCC} we can find the following
inequality in the given condition:
$$
\Vert |u|^{k}u\Vert _{B_{p^{\prime },2}^{s-\delta }}\preceq \Vert u\Vert
_{B_{p,2}^{s}}^{k+1}.\eqno(2.8)
$$%
Since the domain of $s$ is an open set, we may choose $s_{\epsilon }$ and $%
\delta _{\epsilon }$ so that $s<s_{\epsilon }<s_{1}$ and $\delta _{\epsilon
}<\delta ,$ and require them satisfy (2.5) and (2.6). So (2.7) gives the
inequality
\[
\Vert |u|^{k}u\Vert _{B_{p^{\prime },2}^{s_{\epsilon }-\delta _{\epsilon
}}}\preceq \Vert u\Vert _{B_{p,2}^{s_{\epsilon }}}^{k+1}.
\]%
Finally, by Proposition 3, we obtain the desired estimate
\[
\Vert |u|^{k}u\Vert _{B_{p^{\prime },q}^{s-\delta }}\preceq \Vert
|u|^{k}u\Vert _{B_{p^{\prime },2}^{s_{\epsilon }-\delta _{\epsilon
}}}\preceq \Vert u\Vert _{B_{p,2}^{s_{\epsilon }}}^{k+1}\preceq \Vert u\Vert
_{B_{p,q}^{s_{1}}}^{k+1}.
\]

\bigskip

Now, with Lemma 1 we can embed the modulation space into the Besov
space and invoke Lemma 2 to obtain the nonlinear estimates on the Besov
spaces. Then use Proportion 3 to transfer the estimate back to the
modulation spaces. This is the following lemma, which is crucial in this
paper.

\begin{lemma}
(Nonlinear estimate in modulation spaces)\newline
Let $1\leq q<\infty ,2\leq p<\infty ,s\in \mathbb{R}$ $[s-r]\leq k$. Assume that $q\in
\lbrack p^{\prime },p],(1-\frac{2}{p})n<r$, and
$$
\max \{r-n(\frac{1}{q}-\frac{1}{p}),\frac{n}{q^{\prime }}-\frac{r}{k}\}<s<%
\frac{n}{q^{\prime }}-\frac{1}{k}(1-\frac{2}{p})n.\eqno(2.9)
$$%
Then we have
$$
\Vert u^{k+1}\Vert _{M_{p^{\prime },q}^{s-r}}\preceq \Vert u\Vert
_{M_{p,q}^{s}}^{k+1}.\eqno(2.10)
$$
\end{lemma}

\textbf{Proof: }By Lemma 1 we have
$$
\Vert u^{k+1}\Vert _{M_{p^{\prime },q}^{s-r}}\preceq \Vert u^{k+1}\Vert
_{B_{p^{\prime },q}^{s-r+\tau (p^{\prime },q)}}.\eqno(2.11)
$$%
Since $r-n(\frac{1}{q}-\frac{1}{p})<s$, we have $s-r+\tau (p^{\prime },q)>0$%
. Using Lemma 2, we obtain
$$
\Vert u^{k+1}\Vert _{B_{p^{\prime },q}^{s-r+(\frac{1}{q}-\frac{1}{p}%
)n}}\preceq \Vert u\Vert _{B_{p,q}^{s-(\frac{1}{p^{\prime }}-\frac{1}{q}%
)n}}^{k+1}.\eqno(2.12)
$$%
Choose $s_{1}=s+\varepsilon $ in (2.6), then $s$ satisfies
$$
k(\frac{1}{p}-\frac{1}{n}(s-(\frac{1}{p^{\prime }}-\frac{1}{q})n)+\frac{1}{p}%
-\frac{1}{n}(r-n(1-\frac{2}{p})-\varepsilon )=\frac{1}{p^{\prime }}.\eqno%
(2.13)
$$%
By (2.13) we have
$$
s=\frac{n}{q^{\prime }}-\frac{r}{k}+\frac{\varepsilon }{k}\eqno(2.14)
$$%
because $\tau (p^{\prime },q)+\sigma (p,q)=n(1-\frac{2}{p})<r$, where it is
easy to find $0<\varepsilon <r-n(1-\frac{2}{p})$. Combining this with (2.14)
and the condition of Lemma 2, we easily see that the domain of $s$ is
\[
\max \{r-n(\frac{1}{q}-\frac{1}{p}),\frac{n}{q^{\prime }}-\frac{r}{k}\}<s<%
\frac{n}{q^{\prime }}-\frac{n}{k}(1-\frac{2}{p}).
\]%
Finally, we use Lemma 1 again to obtain
\[
\Vert u^{k+1}\Vert _{B_{p^{\prime },q}^{s-r+(\frac{1}{q}-\frac{1}{p}%
)n}}\preceq \Vert u\Vert _{B_{p,q}^{s-(\frac{1}{p^{\prime }}-\frac{1}{q}%
)n}}^{k+1}\preceq \Vert u\Vert _{M_{p,q}^{s}}^{k+1}.
\]%
The lemma is proved.

\begin{remark}
The condition $q\in \lbrack p^{\prime },p]$ is not necessary, but only for
continence in the calculation. So it does not mean that $q$ must be equal to
$2$ when $p=2$. \ In fact, we can find a larger domain of $q$ when $p=2$. \
More precisely, with the same method as above, we may obtain the estimate
\[
\Vert u^{k+1}\Vert _{M_{2,q}^{s-r}}\preceq \Vert u\Vert
_{M_{2,q}^{s}}^{k+1},
\]%
for $q<2\ \ $and$\ $
$$
max\{r+\frac{n}{2}-\frac{n}{q},\frac{n}{2}-\frac{1}{k}[r-(\frac{1}{q}-\frac{1%
}{2})n]\}<s<\frac{n}{2},~~~~~~n(\frac{1}{q}-\frac{1}{2})<r\eqno(2.15)
$$%
or for $q>2$ $\ $and$~$
$$
max\{r,\frac{n}{q^{\prime }}-\frac{1}{k}[r-n(\frac{1}{2}-\frac{1}{q})]\}<s<%
\frac{n}{q^{\prime }},~~~~~~~~~n(\frac{1}{2}-\frac{1}{q})<r.\eqno(2.16)
$$
\end{remark}

\begin{remark}
For $1\leq p<2$, if we switch $p$ and $p^{\prime }$ in the condition of
Lemma 3, the similar conclusion will be obtained, that is
$$
\Vert |u|^{k}u\Vert _{B_{p,q}^{s-\delta }}\preceq \Vert u\Vert
_{B_{p^{\prime },q}^{s_{1}}}^{k+1}.\eqno(2.17)
$$
\end{remark}

Note that Lemma 3 addresses only the estimate for $k\notin \mathbb{Z}$.
When $k\in \mathbb{Z^{+}}$, we can find the following result in \cite{TI} ,
which will be useful in the proof of Theorem 4

\begin{lemma}
Let $s\geq 0,1\leq p,q,p_{i},q_{i}\leq \infty (i=1,2,3,4)$ \ satisfy
$$
\frac{1}{p}=\frac{1}{p_{1}}+\frac{1}{p_{2}}=\frac{1}{p_{3}}+\frac{1}{p_{4}}%
,~~~~~~~~\frac{1}{q}+1=\frac{1}{q_{1}}+\frac{1}{q_{2}}=\frac{1}{q_{3}}+\frac{%
1}{q_{4}},\eqno(2.18)
$$%
We have
$$
\Vert uv\Vert _{M_{p,q}^{s}}\preceq \Vert u\Vert _{M_{p_{1},q_{1}}^{s}}\Vert
v\Vert _{M_{p_{2},q_{2}}}+\Vert u\Vert _{M_{p_{3},q_{3}}}\Vert v\Vert
_{M_{p_{4},q_{4}}^{s}}\eqno(2.19)
$$

This conclusion also holds for $l_{\Box }^{s,q}(L^{r}L^{p})$. That is, for $%
1\leq r,r_{i}\leq \infty $ $(i=1,2,3,4)$ satisfying
\[
\frac{1}{r}=\frac{1}{r_{1}}+\frac{1}{r_{2}}=\frac{1}{r_{3}}+\frac{1}{r_{4}},
\]
we have
$$
\Vert uv\Vert _{l_{\Box }^{s,q}(L^{r}L^{p})}\preceq \Vert u\Vert _{l_{\Box
}^{s,q_{1}}(L^{r_{1}}L^{p_{1}})}\Vert v\Vert _{l_{\Box
}^{0,q_{2}}(L^{r_{2}}L^{p_{2}})}+\Vert u\Vert _{l_{\Box
}^{0,q_{3}}(L^{r_{3}}L^{p_{3}})}\Vert v\Vert _{l_{\Box
}^{s,q_{4}}(L^{r_{4}}L^{p_{4}})}\eqno(2.20)
$$
\end{lemma}

The second task of this paper is to find the global solution of equation
(1.11) with small initial value. As we all know, the crucial part of the
proof is the estimate of nonlinear part. To this end, we also need the
estimate of fractional integral operator in modulation spaces. Recall that
the fractional integral operator is define by
\[
I_{\alpha }(f)(x)=\int_{R^{n}}\frac{f(x-y)}{|y|^{n-\alpha }}dy.
\]

\begin{lemma}
(Boundedness of fractional integral operator in modulation space)\cite{ST}%
\newline
Let $0<\alpha <n$ and $1<p_{1},p_{2},q_{1},q_{2}<\infty $. The fractional
integral operator $I_{\alpha }$ is bounded from $M_{p_{1},q_{1}}^{s}(R^{n})$
to $M_{p_{2},q_{2}}^{s}(R^{n})$ \ or from $l_{\Box
}^{s,q_{1}}(L^{r}L^{p_{1}})$ to $l_{\Box }^{s,q_{2}}(L^{r}L^{p_{2}})$ \ if
and only if
$$
\frac{1}{p_{1}}\leq \frac{1}{p_{2}}-\frac{\alpha }{n}~~\text{and}~~\frac{1}{%
q_{1}}\leq \frac{1}{q_{2}}+\frac{\alpha }{n}.\eqno(2.21)
$$
\end{lemma}

\section{Proof of the main Theorems}

Before we present the proofs, we need to state the Stricharz estimates of
NLKG in modulation spaces. This estimate on the modulation spaces was
obtained in \cite{BHH}.

\begin{lemma}
(Strichart estimate of NLKG in modulation spaces \cite{BHH}).\newline
Let $\ 2\leq p<\infty ,1\leq q<\infty ,\gamma \geq 2\vee (2/\delta )$, where
$\ \alpha $ and $\delta $ are defined in (1.9). We have following the
estimates:
$$
\Vert K^{\prime }(t)f\Vert _{M_{p,q}^{-\alpha }}\preceq (1+t)^{-\delta
}\Vert f\Vert _{M_{p^{\prime },q}},\eqno(3.1)
$$%
$$
\Vert K^{\prime }(t)f\Vert _{l_{\Box }^{-\alpha /2,q}(L^{\gamma
}(R,L^{p}))}\preceq \Vert f\Vert _{M_{2,q}},\eqno(3.2)
$$%
$$
\Vert \int_{0}^{t}K(t-\tau )fd\tau \Vert _{l_{\Box }^{-\alpha
/2,q}(L^{\gamma }(R,L^{p}))}\preceq \Vert f\Vert _{l_{\Box
}^{-1,q}(L^{1}(R,L^{2}))},\eqno(3.3)
$$%
$$
\Vert \int_{0}^{t}K(t-\tau )fd\tau \Vert _{l_{\Box }^{-\alpha
/2,q}(L^{\gamma }(R,L^{p}))}\preceq \Vert f\Vert _{l_{\Box }^{\alpha
/2-1,q}(L^{\gamma ^{\prime }}(R,L^{p^{\prime }}))}.\eqno(3.4)
$$%
In addition, if $q\in \lbrack \gamma ,\gamma ^{\prime }]$, then we have
$$
\Vert K^{\prime }(t)f\Vert _{L^{\gamma }(R,M_{p,q}^{-\alpha /2})}\preceq
\Vert f\Vert _{M_{2,q}},\eqno(3.5)
$$%
$$
\Vert \int_{0}^{t}K(t-\tau )fd\tau \Vert _{L^{\gamma }(R,M_{p,q}^{-\alpha
/2})}\preceq \Vert f\Vert _{L^{1}(R,M_{2,q}^{-1})},\eqno(3.6)
$$%
$$
\Vert \int_{0}^{t}K(t-\tau )fd\tau \Vert _{L^{\gamma }(R,M_{p,q}^{-\alpha
/2})}\preceq \Vert f\Vert _{L^{\gamma ^{\prime }}(R,M_{p,q}^{\alpha /2-1})},%
\eqno(3.7)
$$
\end{lemma}

We also need the following boundedness of \ $K$ \ and \ $K^{\prime }$ \ on
the modulation spaces.

\begin{lemma}
\cite{BKOR}. Let $1\leq p,q<\infty ,s\in \mathbb{R}$ . We have the following
inequalities:
$$
\Vert K^{\prime }(t)f\Vert _{M_{p,q}^{s}}\preceq (1+t)^{n|\frac{1}{2}-\frac{1%
}{p}|}\Vert f\Vert _{M_{p,q}^{s}},\eqno(3.8)
$$%
$$
\Vert K(t)f\Vert _{M_{p,q}^{s}}\preceq (1+t)^{n|\frac{1}{2}-\frac{1}{p}%
|}\Vert f\Vert _{M_{p,q}^{s-1}},\eqno(3.9)
$$%
where
\[
K(t)=\frac{\sin t(I-\Delta )^{\frac{1}{2}}}{(I-\Delta )^{\frac{1}{2}}}%
,~~~~~~K^{\prime }(t)=\cos (I-\Delta )^{\frac{1}{2}}.
\]
\end{lemma}

\textbf{Proof of Theorem 1. }Let\textbf{\ \ }$\delta ,\alpha $ \ and \ $%
\theta $ \ be defined in \ (1.9). Consider the mapping
\[
\Phi :u\rightarrow K^{\prime }(t)u_{0}+K(t)u_{1}-\int_{0}^{t}K(t-\tau
)|u|^{k}ud\tau
\]%
on the Banach space
$$
X_{1}=L^{\infty }(0,T;M_{2,q}^{s})\bigcap L^{\rho }(0,T;M_{2,q}^{s-\beta })%
\eqno(3.10)
$$%
where $\rho =\frac{2}{\delta }$ and $\beta =\frac{\alpha }{2}$. For all $%
2<p\leq \infty $, we can choose $\ \theta $ such that $\delta (p)<1$. So by
Lemma 6, we have
$$
\Vert \Phi \left( u\right) \Vert _{X_{1}}\preceq \Vert u_{1}\Vert
_{M_{p,q}^{s}}+\Vert u_{0}\Vert _{M_{p,q}^{s-1}}+\Vert |u|^{k}u\Vert
_{L^{1}(0,T;M_{2,q}^{s-1})}.\eqno(3.11)
$$%
Choosing $r=1$ in Remark 6 and using the H\"{o}lder inequality, we obtain
that the nonlinear term
$$
\Vert |u|^{k}u\Vert _{L^{1}(0,T;M_{2,q}^{s-1})}\preceq \Vert u\Vert
_{L^{k+1}(0,T;M_{2,q}^{s})}^{k+1}\preceq T\Vert u\Vert _{L^{\infty
}(0,T;M_{2,q}^{s})}^{k+1}.\eqno(3.12)
$$%
Combining (3.11) and (3.12), we have that
$$
\Vert \Phi \left( u\right) \Vert _{X_{1}}\preceq \Vert u_{1}\Vert
_{M_{p,q}^{s}}+\Vert u_{0}\Vert _{M_{p,q}^{s-1}}+T\Vert |u|^{k}u\Vert
_{X_{1}}.\eqno(3.13)
$$%
On the other hand, in a similar estimate it is easy to check
\[
\Vert \Phi \left( u\right) -\Phi \left( v\right) \Vert _{X_{1}}\preceq
T\left( \left\Vert u\right\Vert _{X_{1}}^{k}+\left\Vert v\right\Vert
_{X_{1}}^{k}\right) \left\Vert u-v\right\Vert _{X_{1}}
\]
Denote
\[
B_{M}\text{ }=\{u\in X_{1}:\Vert u\Vert _{X_{1}}\leq M\}.
\]%
We choose $\ M,T>0\ $\ such that $\ $%
\[
\Phi :B_{M}\text{ }\rightarrow B_{M}\text{ }
\]%
is an onto mapping and \
\[
\Vert \Phi \left( u\right) -\Phi \left( v\right) \Vert _{X_{1}}\leq \frac{1}{%
2}\left\Vert u-v\right\Vert _{X_{1}}.
\]%
Thus, we complete the proof of theorem by the standard method of contraction
mapping.

In the proofs of Corollary 1 and Theorem 2, we can not use the Stricharz
estimate in the space
\[
X_{2}=L^{\gamma }(0,T;M_{p,q}^{s}).
\]%
Hence we will invoke the boundedness of Klein-Gordon semigroup in modulation
spaces to estimate the linear part.

In the proof of corollary 1, we first choose $\theta =0$ in (1.9) such
that $\alpha =\delta =0$ in (3.1). So we obtain $\Vert K(t)f\Vert
_{M_{p,q}^{s}}\preceq \Vert f\Vert _{M_{p^{\prime },q}^{s}}$. Now by Lemma
7, we obtain
\[
\Vert \Phi \left( u\right) \Vert _{X_{2}}\preceq (1+T)^{n|\frac{1}{2}-\frac{1%
}{p}|}(\Vert u_{0}\Vert _{M_{p,q}^{s}}+\Vert u_{1}\Vert
_{M_{p,q}^{s-1}})+\Vert \int_{0}^{t}|u|^{k}ud\tau \Vert _{L^{\gamma
}(0,T;M_{p^{\prime },q}^{s-1})}.
\]%
Choosing $r=1$ in Lemma 3 and using H\"{o}lder's inequality, we have
$$
\Vert \Phi \left( u\right) \Vert _{X_{2}}\preceq (1+T)^{n|\frac{1}{2}-\frac{1%
}{p}|}(\Vert u_{0}\Vert _{M_{p,q}^{s}}+\Vert u_{1}\Vert
_{M_{p,q}^{s-1}})+T^{1-\frac{k}{\gamma }}\Vert u\Vert _{X_{2}}^{k+1}.\eqno%
(3.14)
$$

The rest of the proof is the same as that of Theorem 1 with the help of
the method of contraction mapping.

\bigskip

To prove Theorem 2, we use Lemma 7 again to obtain
\begin{eqnarray*}
\Vert \Phi (u)\Vert _{X_{2}} &\preceq &\Vert 1+t^{(\frac{1}{2}-\frac{1}{p}%
)n}\Vert _{L^{\gamma }(0,T)}(\Vert u_{1}\Vert _{M_{p,q}^{s}}+\Vert
u_{0}\Vert _{M_{p,q}^{s-1}})+\Vert \int_{0}^{t}K(t-\tau )|u|^{k}ud\tau \Vert
_{X_{2}} \\
&\preceq &(T^{\frac{1}{\gamma }}+T^{(\frac{1}{2}-\frac{1}{p})n+\frac{1}{%
\gamma }})(\Vert u_{1}\Vert _{M_{p,q}^{s}}+\Vert u_{0}\Vert
_{M_{p,q}^{s-1}})+\Vert \int_{0}^{t}K(t-\tau )|u|^{k}ud\tau \Vert _{X_{2}}.
\end{eqnarray*}%
By choosing $r=1$ in Remark 6 and Remark 7, we use H\"{o}lder's
inequality to obtain the following estimate for the nonlinear term:
\begin{eqnarray*}
\Vert \int_{0}^{t}K(t-\tau )|u|^{k}ud\tau \Vert _{X_{2}} &\preceq &\Vert
\int_{0}^{t}[1+(t-\tau )]^{n(\frac{1}{p}-\frac{1}{2})}\Vert u^{k+1}\Vert
_{M_{p,q}^{s-1}}d\tau \Vert _{L^{\gamma }(0,T)} \\
&\preceq &\Vert \int_{0}^{t}[1+(t-\tau )]^{n(\frac{1}{p}-\frac{1}{2})}\Vert
u\Vert _{M_{p,q}^{s}}^{k+1}d\tau \Vert _{L^{\gamma }(0,T)} \\
&\preceq &\Vert u\Vert _{L^{\gamma }(0,T;M_{p,q}^{s})}^{k+1}\cdot \Vert t^{%
\frac{\gamma -k-1}{\gamma }}(1+t^{1+n(\frac{1}{p}-\frac{1}{2})})\Vert
_{L^{\gamma }(0,T)} \\
&\preceq &(1+T^{1+n(\frac{1}{p}-\frac{1}{2})})T^{\frac{\gamma -k}{\gamma }%
}\Vert u\Vert _{X_{2}}^{k+1}.
\end{eqnarray*}%
If we first assume $T<1$, from the above estimates we obtain that
$$
\Vert \Phi (u)\Vert _{X_{2}}\leq C_{T}[(\Vert u_{1}\Vert
_{M_{p,q}^{s}}+\Vert u_{0}\Vert _{M_{p,q}^{s-1}})]+\Vert u\Vert
_{X_{2}}^{k+1}.\eqno(3.15)
$$%
Also, a similar method gives
$$
\Vert \Phi (u)-\Phi (v)\Vert _{X_{2}}\leq C_{T}(\Vert u\Vert
_{X_{2}}^{k}+\Vert v\Vert _{X_{2}}^{k})\Vert u-v\Vert _{X_{2}},\eqno(3.16)
$$%
where the constant \ $C_{T}$ $\rightarrow 0,$ \ as \ $T\rightarrow 0.$ Now,
by (3.15) and (3.16), the contraction mapping yields the conclusion of
Theorem 2.\newline

\bigskip

\textbf{Proof of Theorem 3. } We denote the space
\[
X=L^{\infty }(R;M_{2,q}^{s})\bigcap L^{k+2}(R;M_{p,q}^{s-\frac{\alpha }{2}%
}),
\]%
where $\alpha $ and $\delta $ are defined in (1.9), and $k+2\geq 2\vee (%
\frac{2}{\delta })$. Since $\delta <1,q\in \lbrack p^{\prime },p]$, we can
choose $q$ such that $(k+2)^{\prime }\leq q\leq \frac{2}{\delta }$. So we
may assume $k+2\geq \frac{2}{\delta }$ for convenience. By Lemma 6, we
have
\[
\Vert \Phi (u)\Vert _{X}\preceq \Vert u_{0}\Vert _{M_{2,q}^{s}}+\Vert
u_{1}\Vert _{M_{2,q}^{s-1}}+\Vert \int_{0}^{t}K(t-\tau )|u|^{k}ud\tau \Vert
_{L^{k+2}(R;M_{p,q}^{s-\frac{\alpha }{2}})}.
\]%
The last term above can be estimated in the following by using Lemma 3 and
choosing $\ r=1-\alpha $, $s=s-\frac{\alpha }{2}$ in the lemma. An easy
computation gives
\[
\Vert \int_{0}^{t}K(t-\tau )|u|^{k}ud\tau \Vert _{L^{k+2}(R;M_{p,q}^{s-\frac{%
\alpha }{2}})}\preceq \Vert u^{k+1}\Vert _{L^{\frac{k+2}{k+1}%
}(R;M_{p^{\prime },q}^{s+\frac{\alpha }{2}-1})}\preceq \Vert u\Vert
_{L^{k+2}(R;M_{p^{\prime },q}^{s-\frac{\alpha }{2}})}^{k+1}.
\]%
So, we have
\[
\Vert \Phi (u)\Vert _{X}\preceq \Vert u_{0}\Vert _{M_{2,q}^{s}}+\Vert
u_{1}\Vert _{M_{2,q}^{s-1}}+\Vert u\Vert _{X}^{k+1}.
\]%
By the standard method used in the proofs of Theorem 1 and Theorem 2, we
can obtain the existence and uniqueness of global solution if the initial
value $\Vert u_{0}\Vert _{M_{2,q}^{s}}+\Vert u_{1}\Vert _{M_{2,q}^{s-1}}$ is
small enough.

Finally, we find the domain of $\ k$. We notice that when we use Lemma 3, $%
\ $the index $p$ should satisfy
\[
n(1-\frac{2}{p})<1-\alpha .
\]%
This gives
\[
2n(\frac{1}{2}-\frac{1}{p})<1-(n+1)\theta (\frac{1}{2}-\frac{1}{p})
\]%
by a simply calculation. So we find that the domain of $p$ is $(\frac{1}{2}-%
\frac{1}{p})>\frac{1}{2n+(n+1)\theta }$.\newline
The condition $k+2\geq \frac{2}{\delta }$ generates
\[
k\geq \frac{2}{\delta }-2>\frac{4}{\theta }+\frac{2}{n}.
\]

\bigskip

\textbf{Proof of Theorem 4. }Checking the above proof, we know that the
global existence theory for small initial data is a straightforward result
of the nonlinear estimate. Thus the main issue is to obtain an estimate for
the nonlinear part. We use $Y$ to denote the space
\[
Y=l_{\Box }^{s,q}(L^{\infty }(R;L^{p}(R^{n})))\bigcap l_{\Box }^{s-\frac{%
\alpha }{2},q}(L^{4}(R;L^{p}(R^{n})).
\]
By Lemma 4, we obtain the following estimate for the nonlinear part:
\begin{eqnarray*}
\Vert (|x|^{-\mu }\ast |u|^{2})u\Vert _{l_{\Box }^{s+\frac{\alpha }{2}%
-1,q}(L_{T}^{4/3}L^{p^{\prime }})} &\preceq &\Vert |x|^{-\mu }\ast
|u|^{2}\Vert _{l_{\Box }^{s+\frac{\alpha }{2}-1,q_{1}}(L_{T}^{2}L^{p_{1}})}%
\Vert u\Vert _{l_{\Box }^{0,q_{2}}(L_{T}^{4}L^{p_{2}})} \\
&+&\Vert |x|^{-\mu }\ast |u|^{2}\Vert _{l_{\Box
}^{0,q_{3}}(L_{T}^{4}L^{p_{3}})}\Vert u\Vert _{l_{\Box }^{s+\frac{\alpha }{2}%
-1,q_{4}}(L_{T}^{2}L^{p_{4}})}=I+II.
\end{eqnarray*}%
We will only estimate term $I$, since the second term $II$ can be estimated
in the same way. Using Proposition 2, Lemma 4 and Lemma 5, we have
\begin{eqnarray*}
I &\preceq &\Vert u^{2}\Vert _{l_{\Box }^{s+\frac{\alpha }{2}%
-1,q_{5}}(L_{T}^{2}L^{p_{5}})}\Vert u\Vert _{l_{\Box
}^{0,q_{2}}(L_{T}^{4}L^{p_{2}})} \\
&\preceq &\Vert u\Vert _{l_{\Box }^{s+\frac{\alpha }{2}%
-1,q_{6}}(L_{T}^{4}L^{p_{6}})}\Vert u\Vert _{l_{\Box
}^{0,q_{7}}(L_{T}^{4}L^{p_{7}})}\Vert u\Vert _{l_{\Box
}^{0,q_{2}}(L_{T}^{4}L^{p_{2}})} \\
&\preceq &\Vert u\Vert _{l_{\Box }^{s-\frac{\alpha }{2}%
,q}(L_{T}^{4}L^{p})}^{3}.
\end{eqnarray*}%
Similarly, we can obtain
\[
II\preceq \Vert u\Vert _{l_{\Box }^{s-\frac{\alpha }{2}%
,q}(L_{T}^{4}L^{p})}^{3}.
\]%
Therefore, we have
\[
\Vert (|x|^{-\mu }\ast |u|^{2}u)\Vert _{l_{\Box }^{s+\frac{\alpha }{2}%
-1,q}(L_{T}^{4/3}L^{p^{\prime }})}\preceq \Vert u\Vert _{l_{\Box }^{s-\frac{%
\alpha }{2},q}(L_{T}^{4}L^{p})}^{3}\preceq \Vert u\Vert _{Y}^{3}.
\]%
Again, by the standard method of contraction mapping, we prove the
conclusion of the theorem.

Finally we check the range of \ $\mu .$ In the above proof, we notice that
the conditions in the estimate of \ $I$ \ (the same condition in the
estimate of \ $II)$ imply that $p$ should satisfy
\[
\frac{1}{p_{1}}+\frac{1}{p_{2}}=\frac{1}{p^{\prime }},
\]%
\[
\frac{1}{p_{6}}+\frac{1}{p_{7}}=\frac{1}{p_{5}},
\]%
\[
\frac{1}{p_{1}}\leq \frac{1}{p_{5}}-\frac{n-\mu }{n},
\]%
and
\[
p_{2},p_{6},p_{7}\leq p.
\]%
These clearly yield that
\[
2n(1-\frac{2}{p})\leq \mu .
\]%
Also, we note that $\ q$ \ should satisfy
\[
\frac{1}{q_{1}}+\frac{1}{q_{2}}=\frac{1}{q}+1,
\]%
\[
\frac{1}{q_{6}}+\frac{1}{q_{7}}=\frac{1}{q_{5}}+1,
\]%
\[
\frac{1}{q_{1}}\leq \frac{1}{q_{5}}+\frac{n-\mu }{n},
\]%
\[
\frac{n}{q_{7}},\frac{n}{q_{2}}<s-\frac{\alpha }{2}+\frac{n}{q},
\]%
and
\[
s+\frac{\alpha }{2}-1+\frac{n}{q_{6}}<s-\frac{\alpha }{2}+\frac{n}{q}.
\]%
Thus, a direction computation gives that
\[
\mu \leq 2s+\frac{2n}{q}+1-2\alpha -n.
\]%
So, the domain of $\mu $ \ is
\[
2n(1-\frac{2}{p})\leq \mu \leq 2s+\frac{2n}{q}+1-2\alpha -n.
\]

To compare the low bound of \ $\mu $ \ to that in \cite{MZ} when \ $n=3,$ we
choose $\theta =1$ in (1.9). The condition $4\geq \frac{2}{\delta }$ means
that $4(1-\frac{2}{p})\geq \frac{1}{2}$. Since $n=3$, the value $2n(1-\frac{2%
}{p})$ \ is at least $\frac{3}{4}$.

\end{document}